\documentclass{mcom-l}

\usepackage{amssymb}
\usepackage{enumitem}
\usepackage{booktabs}
\usepackage{hyperref}

\newtheorem{theorem}{Theorem}[section]
\newtheorem{lemma}[theorem]{Lemma}
\newtheorem{proposition}[theorem]{Proposition}
\newtheorem{corollary}[theorem]{Corollary}
\newtheorem*{thmA}{Theorem A}
\newtheorem*{thmB}{Theorem B}

\theoremstyle{definition}
\newtheorem{definition}[theorem]{Definition}

\theoremstyle{remark}
\newtheorem{remark}[theorem]{Remark}

\numberwithin{equation}{section}

\newcommand{\Z}{\mathbb{Z}}
\newcommand{\Q}{\mathbb{Q}}
\newcommand{\R}{\mathbb{R}}
\newcommand{\one}{\mathbf{1}}
\newcommand{\covol}{\operatorname{covol}}
\newcommand{\dist}{\operatorname{dist}}
\newcommand{\spn}{\operatorname{span}}
\newcommand{\nZ}[1]{\lVert #1\rVert_{\Z}}

\begin{document}

\title[View-obstruction for spheres in dimensions 6 and 7]{The view-obstruction problem for spheres\\ in dimensions six and seven}

\author{Vineet Kahlon}
\address{Google, New York, NY, USA}



\begin{abstract}
Place closed Euclidean balls of radius $r$ at the points of $\tfrac12\one+\Z^n$. The view-obstruction
problem for spheres asks for the least $r$ such that every ray from the origin into the open positive
orthant meets a ball. Equivalently, it asks for $\Delta_n=\sup_v\min_{t\in\R}\sum_i\nZ{tv_i-\tfrac12}^2$
over positive integer vectors $v$. The values for $n\le5$ are classical. We prove
$\Delta_6=107/312$, attained only by permutations of $(12,9,7,5,3,2)$, and $\Delta_7=73/180$, attained only by permutations of
$(11,10,8,7,4,3,1)$. The proof reduces the supremum to a finite search. The key facts are an exact
lower bound $d_n(v)\ge\delta_U^2$ for every rational subspace $U\perp v$, a matching fiber upper bound, and
Banaszczyk's transference theorem, applied along the Korkine--Zolotarev flag of the dual lattice $v^\perp\cap\Z^n$.
The search is quotiented by signed coordinate permutations and carried out entirely in exact rational
arithmetic. It terminates with a certificate: no node of the search is left unresolved, and every
surviving candidate is evaluated exactly. The method also recovers the known values for $n=3,4,5$.
\end{abstract}

\maketitle

\section{Introduction}

\subsection{The problem}

Imagine standing at the origin of $\R^n$ in a regular forest: a tree of radius $r$ grows at every point of
the shifted lattice $c+\Z^n$, where $c=\tfrac12\one=(\tfrac12,\dots,\tfrac12)$. Look out into the open
positive orthant. When is the view blocked in every direction? For tiny $r$ some rays escape to infinity;
for large $r$ every ray eventually hits a tree. The \emph{view-obstruction constant} is the threshold
radius, and this paper determines it for $n=6$ and $n=7$.

The shift by $c$ matters. Trees at the integer points would be hit at once, since every ray starts at a
lattice point. The half-integer points are the positions farthest from the lattice, so they give the most
interesting obstruction problem. With cubes in place of balls, the same question is the well-known lonely
runner problem (Section~\ref{subsec:lrc}).

Two simple observations turn the problem into number theory. First, a ray whose direction is not proportional
to an integer vector winds densely around a subtorus of $\R^n/\Z^n$. It therefore comes at least as close to the trees as the rational rays inside that subtorus do, so it is
never harder to block (Proposition~\ref{prop:reduction}). Second, along the ray through an integer vector $v$ the configuration
repeats with period one. So only a single period needs to be examined. With
$\nZ{x}=\min_{k\in\Z}|x-k|$ we are led to
\begin{equation}\label{eq:d}
d_n(v)=\inf_{t\in\R}\sum_{i=1}^n\nZ{tv_i-\tfrac12}^2=\inf_{t\in\R}\dist(tv,\,c+\Z^n)^2,
\qquad
\Delta_n=\sup_{v\in\Z^n_{>0}}d_n(v).
\end{equation}
Here $d_n(v)$ is the squared radius needed to block the single ray through $v$, and $\sqrt{\Delta_n}$ is the
view-obstruction constant (Proposition~\ref{prop:reduction}). In the normalization of
\cite{DHW94,DHW96}, which uses balls of diameter~$1$ dilated by $\nu$, one has $\nu(n)^2=4\Delta_n$.

A first estimate shows how much room there is. For each $i$, the fractional part of $tv_i$ is uniformly
distributed as $t$ runs over $[0,1]$. So the average of $\sum_i\nZ{tv_i-\frac12}^2$ over one period is exactly
$n/12$, and hence
\begin{equation}\label{eq:trivial}
\Delta_n\le\frac{n}{12}.
\end{equation}
The difficulty lies in identifying the directions that stay far from the trees for the whole period, and in
proving that no direction does better. The second part is hard because there are infinitely many directions. A long vector $v$ gives a closed
orbit that winds many times around the torus before it repeats, and no finite list of directions is
obviously enough.

\subsection{Cubes and the lonely runner problem}\label{subsec:lrc}

Replace each ball by the cube $c+z+\alpha[-\frac12,\frac12]^n$, of side $\alpha$. The ray through $v$ meets
such a cube when, for some $t$, every coordinate satisfies $\nZ{tv_i-\frac12}\le\alpha/2$. Equivalently,
$\nZ{tv_i}\ge\frac{1-\alpha}2$ for all $i$. Read $t$ as time and $tv_i \bmod 1$ as the position of a runner
with speed $v_i$ on a circular track of length one. The ray is blocked exactly when, at some moment, every
runner is at distance at least $\frac{1-\alpha}2$ from the start. The \emph{lonely runner
conjecture} of Wills~\cite{Wil67} and Cusick~\cite{Cus73} asserts that for every $v\in\Z^n_{>0}$ there is a
time with $\nZ{tv_i}\ge\frac1{n+1}$ for all $i$. The bound is attained for $v=(1,2,\dots,n)$. In the
language above, the conjecture says that the view-obstruction constant for cubes is $(n-1)/(n+1)$. In the
usual formulation this is the statement for $n+1$ runners, one of them standing still.

The conjecture is known for up to seven runners \cite{BHK01,BS08}. Recently, computer-assisted proofs
have extended it to eight runners \cite{Ros25a}, and then to nine and ten runners
\cite{Ros25b,Tra25}. So the cube problem is settled for $n\le9$. For connections with covering radii of
lattice zonotopes, see Henze and Malikiosis~\cite{HM16}.

Cubes and balls differ in an essential way. For cubes the obstruction condition is coordinatewise: each
runner only has to be in a fixed arc, independently of the others. This makes the problem combinatorial,
a question about covering a period by unions of intervals. For balls the coordinates interact through the
sum of squares. The quantity $d_n(v)$ is a closest-vector distance in a projected lattice, and the natural
tools come from the geometry of numbers: duality, transference and reduction theory. The extremal
directions reflect this difference. For cubes the arithmetic progression $(1,2,\dots,n)$ is extremal. For
balls the extremal directions, such as $(12,9,7,5,3,2)$, have no such simple pattern. The recent lonely
runner results and the present paper share a strategy: a theoretical reduction to a finite search,
followed by an exhaustive computer verification. The reductions themselves are quite different.

\subsection{What was known}

The problem goes back to Cusick~\cite{Cus73}. For spheres it was settled by increasingly intricate case
analyses: $n=3$ by Dumir and Hans-Gill~\cite{DH86}, $n=4$ by Dumir, Hans-Gill and Wilker~\cite{DHW94},
and $n=5$ by the same authors~\cite{DHW96}. Each dimension required its own geometric and Diophantine
argument.
As a check on both the method and the literature, the search developed here re-derives the
values for $n=3,4,5$ independently of these case analyses (Section~\ref{subsec:small}).
Together with our results, the values are listed in Table~\ref{tab:values}.

\begin{table}[ht]
\caption{Known values of $\Delta_n$ and extremal directions.}\label{tab:values}
\centering\small
\begin{tabular}{clcc}
\toprule
$n$ & extremal direction $v$ & $\Delta_n$ & $\Delta_n/(n/12)$\\
\midrule
2 & $(2,1)$ & $1/20$ & $0.300$\\
3 & $(3,2,1)$ & $3/28$ & $0.429$\\
4 & $(3,2,1,1)$ & $11/60$ & $0.550$\\
5 & $(9,7,5,3,2)$ & $41/168$ & $0.586$\\
6 & $(12,9,7,5,3,2)$ & $107/312$ & $0.686$\\
7 & $(11,10,8,7,4,3,1)$ & $73/180$ & $0.695$\\
\bottomrule
\end{tabular}
\end{table}

The last column compares $\Delta_n$ with the averaging bound~\eqref{eq:trivial}. The extremal directions are
short vectors with small, mostly distinct entries. Heuristically, the rays that stay farthest from the
trees are those whose coordinates pass through the half-integers at well-separated times within a
period.

\subsection{Results}

\begin{thmA}
$\Delta_6=\dfrac{107}{312}$. A positive primitive $v\in\Z^6$ satisfies $d_6(v)=\Delta_6$ if and only if $v$ is a
permutation of $(12,9,7,5,3,2)$.
\end{thmA}

\begin{thmB}
$\Delta_7=\dfrac{73}{180}$. A positive primitive $v\in\Z^7$ satisfies $d_7(v)=\Delta_7$ if and only if $v$ is a
permutation of $(11,10,8,7,4,3,1)$.
\end{thmB}

The lower bounds in both theorems come from exhibiting one direction: the time $t=19/312$, respectively
$t=11/180$, already witnesses the value (Section~\ref{sec:reduction}). The content of the theorems is the
upper bound. To prove it, we reduce ``no direction beats $T$'' to a finite, exactly checkable computation.
The same reduction, run at the classical thresholds, recovers the known values for $n=3,4,5$.

\subsection{The idea of the proof}

Our approach replaces dimension-specific case analysis with one uniform mechanism. It rests on three ideas.

\emph{1. A good direction forces a short orthogonal vector.} Project $\Z^n$ orthogonally along $v$. The
value $d_n(v)$ is then the squared distance from the projected point $c$ to a lattice $\Gamma_v$ of rank
$n-1$. If $d_n(v)$ is large, $\Gamma_v$ has a large covering radius. By a transference theorem of Banaszczyk,
its dual lattice $L_v=v^\perp\cap\Z^n$ then contains a short vector. So every direction with $d_n(v)>T$ is
orthogonal to one of finitely many short integer vectors, which can be listed in advance.

\emph{2. Repeat inside the orthogonal complement, and account for what has been paid.} Once we know short
vectors $w_1,\dots,w_k$ orthogonal to $v$, spanning $U$, part of the distance is already determined by $U$
alone. This is the \emph{offset} $\delta_U^2$, the squared distance from $c$ to the lattice $\Z^n$ projected
onto $U$. It is a lower bound for $d_n(v)$ (Lemma~\ref{lem:lower}). The remaining gap $T-\delta_U^2$ must come
from the complementary directions, and transference again forces a short vector there. After $n-2$ steps
only a rank-two lattice is left. It contains $v$, and the identity $\|v\|=\det L_v$ bounds $\|v\|$
(Proposition~\ref{prop:leaf}). Every step depends only on $U$, never on the unknown $v$, so the procedure is a
finite search tree.

\emph{3. The only way the recursion can stall is by finding good directions.} The step bound needs a
positive gap $T-\delta_U^2$. If some node has $\delta_U^2\ge T$, the lower bound of idea~2 shows that
\emph{every} direction orthogonal to $U$ already has $d_n\ge T$. Such a node is therefore not a gap in the
argument but evidence against the claimed value (Proposition~\ref{prop:witness}). We record such nodes as
\emph{unresolved}. The resulting dichotomy is clean. Either the search closes with no unresolved node, and
every surviving candidate is checked exactly, which proves $\Delta_n\le T$ (Corollary~\ref{cor:cert}). Or it
produces explicit directions with $d_n\ge T$.

To keep the tree small, we follow one canonical chain of short vectors for each $v$, a Korkine--Zolotarev
(KZ) flag of $L_v$. This comes with strong necessary conditions, such as Hermite's inequality in rank two, and
we quotient by the signed coordinate permutations, which preserve $d_n$. All arithmetic is exact, with
integers and rationals throughout.

The method has no built-in dimensional barrier: the search is finite in every dimension, and only its cost
grows. For $n=7$ the certified run took about $27$ hours on eight processes.

\subsection{Organization}

Section~\ref{sec:reduction} gives the reduction to integer directions, an exact formula for $d_n(v)$, and the
symmetries. Section~\ref{sec:lattices} develops the offset and transference bounds behind ideas~1 and~2.
Section~\ref{sec:kz} introduces KZ flags and proves the completeness theorem and the certificate criterion
(idea~3). Section~\ref{sec:algorithm} turns this into an algorithm and explains why it is exact.
Section~\ref{sec:computations} reports the computations and proves Theorems~A and~B.
Section~\ref{sec:remarks} collects further remarks.

\section{Reduction, exact evaluation and symmetry}\label{sec:reduction}

This section justifies the formulation~\eqref{eq:d}. It also shows that $d_n(v)$ can be computed exactly for
any single $v$, and records the symmetries that the search will exploit.

\begin{proposition}[Reduction to rational directions]\label{prop:reduction}
For $x\in\R^n_{>0}$ let $\rho(x)=\inf_{t>0}\dist(tx,c+\Z^n)$, and let $\kappa_n$ be the infimum of all
$r>0$ such that every ray $\{tx:t>0\}$ with $x\in\R^n_{>0}$ meets a closed ball of radius $r$ centred
in $c+\Z^n$. Then $\kappa_n^2=\sup_x\rho(x)^2=\Delta_n$.
\end{proposition}

\begin{proof}
If $r>\sup\rho$, every ray comes within distance $<r$ of a centre. If $r<\rho(x)$, the ray through $x$
misses every ball. Hence $\kappa_n=\sup_x\rho(x)$.

Let $v\in\Z^n_{>0}$. The function $f(t)=\dist(tv,c+\Z^n)^2$ is $1$-periodic. It is even, because
$-c+\Z^n=c+\Z^n$. It is also continuous. So $\rho(v)^2=\min_{t\in\R}f(t)=d_n(v)$.

For general $x$, let $W$ be the smallest subspace of $\R^n$ that is defined over $\Q$ and contains $x$. By
Kronecker's theorem, the closure in $\R^n/\Z^n$ of $\{tx:t\ge0\}$ is the subtorus $(W+\Z^n)/\Z^n$; the
closure of a half-orbit of a flow on a compact abelian group is a closed subsemigroup, hence a subgroup.
Since the distance is continuous,
\[
\rho(x)=\dist(c+\Z^n,\,W+\Z^n)\le\dist(c+\Z^n,\,\R v+\Z^n)=\sqrt{d_n(v)}
\]
for every $v\in W$. The set $W\cap\R^n_{>0}$ is a nonempty open subset of the rational space $W$, so it
contains a primitive integer vector $v$. Thus $\rho(x)^2\le\Delta_n$ for all $x$, with equality for
suitable integer $x$.
\end{proof}

\begin{lemma}[Exact evaluation]\label{lem:exact}
Let $v\in\Z^n_{>0}$, let $0=t_0<t_1<\dots<t_N=1$ be the distinct points $k/v_i$
($1\le i\le n$, $0\le k\le v_i$), and for $0\le j<N$ put $c^{(j)}_i=\lfloor\tau_jv_i\rfloor+\tfrac12$
with $\tau_j=\tfrac12(t_j+t_{j+1})$. Then
\[
d_n(v)=\min_{0\le j<N}\Big(\|c^{(j)}\|^2-\frac{\langle v,c^{(j)}\rangle^2}{\|v\|^2}\Big).
\]
\end{lemma}

\begin{proof}
For $t\in(t_j,t_{j+1})$ and each $i$, $\lfloor tv_i\rfloor=\lfloor\tau_jv_i\rfloor$. Hence
$\nZ{tv_i-\frac12}=|tv_i-c^{(j)}_i|$, and $f(t)=\|tv-c^{(j)}\|^2$ on that interval. The $j$-th bracket equals
$\min_{s\in\R}\|sv-c^{(j)}\|^2$, and it is $\ge\inf f$ because $\|sv-c'\|^2\ge f(s)$ for every
$c'\in c+\Z^n$. Conversely, let $t^\ast\in[0,1]$ minimize $f$, and take a cell $j$ whose closure contains
$t^\ast$. By continuity, the $j$-th bracket is $\le\|t^\ast v-c^{(j)}\|^2=f(t^\ast)$. (The unconstrained
minimizer $s$ of a cell may lie outside that cell; this does not affect the argument.)
\end{proof}

For instance, $t=19/312$ with $c^{(0)}=c$ gives $d_6(12,9,7,5,3,2)\le 107/312$, and $t=11/180$ gives
$d_7(11,10,8,7,4,3,1)\le 73/180$. Lemma~\ref{lem:exact} shows that equality holds in both cases.

\begin{lemma}[Symmetry and monotonicity]\label{lem:sym}
\begin{enumerate}[label=(\alph*)]
\item Extend $d_n$ by \eqref{eq:d} to all $v\in\Z^n$ with no zero coordinate. For every signed
permutation matrix $S$, $d_n(Sv)=d_n(v)=d_n(|v|)$.
\item $d_{n+1}(v,a)\ge d_n(v)$ for all $a\in\Z_{>0}$. Hence $\Delta_{n+1}\ge\Delta_n$.
\end{enumerate}
\end{lemma}

\begin{proof}
(a) holds because $\nZ{-x-\frac12}=\nZ{x+\frac12}=\nZ{x-\frac12}$ and permutations only reorder the sum.
(b) The extra summand $\nZ{ta-\frac12}^2$ is nonnegative. The vector $(v,a)$ is primitive whenever $v$ is.
\end{proof}

\begin{remark}[Averaging bound]\label{rem:avg}
For $v\in\Z^n_{>0}$ and each $i$, the map $t\mapsto tv_i \bmod 1$ pushes Lebesgue measure on $[0,1]$
forward to Lebesgue measure. Hence $\int_0^1\nZ{tv_i-\frac12}^2\,dt=\int_0^1\nZ{u}^2\,du=\frac1{12}$, and
$d_n(v)\le n/12$. This proves~\eqref{eq:trivial}.
\end{remark}

\section{Lattices, offsets and transference}\label{sec:lattices}

Here we make ideas~1 and~2 of the introduction precise. The guiding picture is this. Projecting along $v$
turns $d_n(v)$ into a closest-vector distance in the lattice $\Gamma_v$. A rational subspace $U\perp v$ splits
this distance into a part that is visible inside $U$, the offset, and a remainder in the complement. The
offset bounds $d_n(v)$ from below. The remainder is bounded above by a covering radius, which transference
converts into the existence of a short dual vector.

Throughout, $v\in\Z^n_{>0}$ is primitive. Put
\[
L=L_v=v^\perp\cap\Z^n,\qquad\Gamma=\Gamma_v=P_{v^\perp}(\Z^n),
\]
where $P_X$ denotes orthogonal projection onto a subspace $X$. We use the following standard facts
\cite{Cas59,Mar03}.
\begin{enumerate}[label=(F\arabic*)]
\item\label{F1} $L$ is a lattice of rank $n-1$ with $\det L=\|v\|$, and $\Gamma$ is a lattice in $v^\perp$
with dual $\Gamma^\ast=L$.
\item\label{F2} Let $\Lambda$ be a lattice spanning $V$, and let $W\subseteq V$ be spanned by vectors of
$\Lambda^\ast$. Then $S=V\cap W^\perp$ is spanned by vectors of $\Lambda$, and
$(\Lambda\cap S)^\ast=P_S(\Lambda^\ast)$.
\item\label{F3} Let $M$ be a lattice and $W$ a subspace spanned by vectors of $M$. Then
$\det M=\det(M\cap W)\cdot\det P_{W^\perp}(M)$.
\item\label{F4} (Banaszczyk \cite{Ban93}) A lattice $\Lambda$ of rank $m$ with covering radius $\mu(\Lambda)$
satisfies $\mu(\Lambda)\,\lambda_1(\Lambda^\ast)\le m/2$. For $m=1$ this is an equality.
\item\label{F5} (Hermite, $\gamma_2=2/\sqrt3$) A rank-$2$ lattice with minimum $\lambda$ and Gram--Schmidt
norms $\lambda,h$ satisfies $h\ge\frac34\lambda$.
\end{enumerate}

\begin{definition}
For a linear subspace $U\subseteq\R^n$, the \emph{offset} of $U$ is
$\delta_U^2=\min_{z\in\Z^n}\|P_U(c-z)\|^2$. If $U$ is spanned by integer vectors, then
$P_U(\Z^n)=(\Z^n\cap U)^\ast$ and $\delta_U^2$ is the squared distance from $P_Uc$ to this lattice. It depends
only on $U$.
\end{definition}

The next lemma is elementary but central. Whatever $v$ is, the part of $c-z$ visible inside $U$ is fixed
by $U$ alone. So the offset is a floor for the whole distance.

\begin{lemma}[Offset lower bound]\label{lem:lower}
If $U\subseteq v^\perp$, then $d_n(v)\ge\delta_U^2$.
\end{lemma}

\begin{proof}
By \eqref{eq:d}, $d_n(v)=\min_{z}\|P_{v^\perp}(c-z)\|^2$. For every $z$,
$\|P_{v^\perp}(c-z)\|^2=\|P_U(c-z)\|^2+\|P_{v^\perp\cap U^\perp}(c-z)\|^2\ge\delta_U^2$.
\end{proof}

Conversely, if we fix the part of $z$ that realizes the offset and vary $z$ only along $U^\perp$, the
remaining distance is controlled by the covering radius of a smaller lattice.

\begin{lemma}[Fiber upper bound]\label{lem:fiber}
Let $U\subseteq v^\perp$ be spanned by vectors of $L$, with $k=\dim U\le n-2$. Put $S=v^\perp\cap U^\perp$
and $\Lambda_U=P_S(\Z^n\cap U^\perp)$. Then:
\begin{enumerate}[label=(\alph*)]
\item $\Lambda_U=\Gamma\cap S$ is a lattice of rank $n-1-k$ in $S$, and its dual lattice is $P_S(L)=P_{U^\perp}(L)$;
\item $d_n(v)\le\delta_U^2+\mu(\Lambda_U)^2$.
\end{enumerate}
\end{lemma}

\begin{proof}
(a) Since $U\subseteq v^\perp$, we have $P_UP_{v^\perp}=P_U$. So $P_{v^\perp}z\in U^\perp$ if and only if
$z\in U^\perp$, and then $P_{v^\perp}z=P_Sz$. This gives $\Gamma\cap S=\Lambda_U$. The duality statement is
\ref{F2} with $\Lambda=\Gamma$ and \ref{F1}. For $w\in L$, $P_Sw=P_{U^\perp}w$.

(b) Choose $z_0$ with $\|P_U(c-z_0)\|^2=\delta_U^2$. For $y\in\Z^n\cap U^\perp$ we have
$P_U(c-z_0-y)=P_U(c-z_0)$ and $P_S(c-z_0-y)=P_S(c-z_0)-P_Sy$. Minimizing over $y$,
$d_n(v)\le\delta_U^2+\dist(P_S(c-z_0),\Lambda_U)^2\le\delta_U^2+\mu(\Lambda_U)^2$.
\end{proof}

Combining the two bounds with transference gives the basic step. If the gap left after $U$ is positive,
the next projected dual lattice must contain a short vector.

\begin{proposition}[Step bound]\label{prop:step}
In the setting of Lemma~\ref{lem:fiber}, if $d_n(v)>T>\delta_U^2$, then
\[
\lambda_1\big(P_{U^\perp}(L)\big)^2<\frac{(n-1-k)^2}{4\,(T-\delta_U^2)}.
\]
\end{proposition}

\begin{proof}
Lemma~\ref{lem:fiber}(b) gives $\mu(\Lambda_U)^2>T-\delta_U^2$. Now apply \ref{F4} with $m=n-1-k$ and
Lemma~\ref{lem:fiber}(a).
\end{proof}

For $U=0$ this gives $\lambda_1(L)^2<(n-1)^2/(4T)$. Because $\|w\|^2\in\Z$, it follows that
$\lambda_1(L)^2\le\lceil(n-1)^2/(4T)\rceil-1$, which is $18$ for $(n,T)=(6,\frac{107}{312})$ and $22$
for $(7,\frac{73}{180})$. When $\delta_U^2\ge T$, Lemma~\ref{lem:lower} shows that no such bound can hold
in general; this is the origin of the unresolved nodes below.

The step bound concerns the unknown lattice $L_v$. The next lemma shows that the search does not need to
know $v$. It can enumerate the finitely many candidate extensions in the fixed lattice $\Z^n$, and a
genuine extension for $L_v$ will be among them.

\begin{lemma}[Independence from $v$]\label{lem:vfree}
Let $U\subseteq v^\perp$, and let $w_1,\dots,w_k$ be a basis of $\Z^n\cap U$. Then:
\begin{enumerate}[label=(\alph*)]
\item $L\cap U=\Z^n\cap U$, so $\covol(w_1,\dots,w_k)=\det(L\cap U)$;
\item $P_{U^\perp}(L)$ is a saturated sublattice of $P_{U^\perp}(\Z^n)\cong\Z^n/(\Z^n\cap U)$;
\item if $P_{U^\perp}w$ is primitive in $P_{U^\perp}(L)$, then $w_1,\dots,w_k,w$ is a basis of
$\Z^n\cap(U+\R w)$.
\end{enumerate}
\end{lemma}

\begin{proof}
(a) holds because $\Z^n\cap U\subseteq v^\perp$. For (b), $\Z^n/L$ is torsion-free and
$\big(\Z^n/(\Z^n\cap U)\big)/\big(L/(\Z^n\cap U)\big)\cong\Z^n/L$. Part (c) is the basis-extension criterion
applied in $\Z^n/(\Z^n\cap U)$, using (b).
\end{proof}

After $n-2$ steps only a single dimension of $L_v$ remains, and $v$ itself can be bounded.

\begin{proposition}[Leaf bound]\label{prop:leaf}
Let $U\subseteq v^\perp$ with $\dim U=n-2$, let $w_1,\dots,w_{n-2}$ be a basis of $\Z^n\cap U$, and suppose
$d_n(v)>T>\delta_U^2$. Then
\[
\|v\|^2<\frac{\covol(w_1,\dots,w_{n-2})^2}{4\,(T-\delta_U^2)}.
\]
\end{proposition}

\begin{proof}
By \ref{F1}, \ref{F3} and Lemma~\ref{lem:vfree}(a),
$\|v\|=\det L=\covol(w_1,\dots,w_{n-2})\,\lambda_1(P_SL)$, since $P_SL$ has rank one. By
Lemma~\ref{lem:fiber} and the equality case of \ref{F4},
$\lambda_1(P_SL)=1/(2\mu(\Lambda_U))<1/(2\sqrt{T-\delta_U^2})$.
\end{proof}

\section{Korkine--Zolotarev flags and the completeness theorem}\label{sec:kz}

The bounds of Section~\ref{sec:lattices} hold for \emph{any} chain of subspaces built from short vectors. A
search over all such chains would be enormous. We therefore fix one canonical chain for each $v$: at every
step take a shortest vector of the current projected lattice. This is a Korkine--Zolotarev flag. It has two
advantages. The transference bound applies to exactly the vector we choose. And the chain must satisfy
additional necessary conditions, which prune the search heavily.

\begin{definition}
A \emph{KZ flag} of $L$ is a basis $w_1,\dots,w_{n-1}$ of $L$ such that, with
$U_j=\spn(w_1,\dots,w_j)$ and $U_0=0$, the vector $P_{U_{j-1}^\perp}w_j$ is a shortest nonzero vector
of $P_{U_{j-1}^\perp}(L)$ for every $j$. Put $h_j=\|P_{U_{j-1}^\perp}w_j\|^2$. KZ flags exist; see
\cite{LLS90}.
\end{definition}

\begin{lemma}[Necessary conditions]\label{lem:kz}
Let $(w_j)$ be a KZ flag of $L$ and $1\le k\le n-2$.
\begin{enumerate}[label=(\alph*)]
\item $h_{k+1}\ge\tfrac34h_k$.
\item For each $j\le k+1$, $P_{U_{j-1}^\perp}w_j$ is a shortest nonzero vector of
$P_{U_{j-1}^\perp}(\Z^n\cap U_{k+1})$. Equivalently, for $m\in\Z\setminus\{0\}$, no vector of
$m\,w_{k+1}+(\Z^n\cap U_k)$ has projection to $U_{j-1}^\perp$ of squared norm $<h_j$. Only
$|m|<\sqrt{h_j/h_{k+1}}$ needs to be tested.
\end{enumerate}
\end{lemma}

\begin{proof}
(a) Apply \ref{F5} to the lattice generated by $P_{U_{k-1}^\perp}w_k$ and $P_{U_{k-1}^\perp}w_{k+1}$. Its
minimum is $h_k$, because it is a sublattice of $P_{U_{k-1}^\perp}(L)$ containing $P_{U_{k-1}^\perp}w_k$.
(b) $\Z^n\cap U_{k+1}\subseteq L$ by Lemma~\ref{lem:vfree}(a), and a shortest vector of a lattice remains
shortest in any sublattice that contains it. A vector of $m\,w_{k+1}+(\Z^n\cap U_k)$ has component of
squared norm $m^2h_{k+1}$ orthogonal to $U_k\supseteq U_{j-1}$.
\end{proof}

The completeness theorem follows the KZ flag of a hypothetical counterexample step by step, and it
records the constraints the flag must satisfy. It stops either at a leaf, where $v$ is bounded, or at the
first subspace whose offset already reaches $T$.

\begin{theorem}[Completeness]\label{thm:complete}
Let $T>0$, let $v\in\Z^n_{>0}$ be primitive with $d_n(v)>T$, and let $(w_j)$ be a KZ flag of $L_v$. Let
$k^\ast$ be the least $k\in\{1,\dots,n-2\}$ with $\delta_{U_k}^2\ge T$, or $k^\ast=\infty$ if there is none. Then:
\begin{enumerate}[label=(\roman*)]
\item $w_1$ is primitive, has at least two nonzero entries of opposite sign, and satisfies
$\|w_1\|^2\le\lceil(n-1)^2/(4T)\rceil-1$;
\item for $1\le k<\min(k^\ast,n-2)$, the class of $w_{k+1}$ in $\Z^n/(\Z^n\cap U_k)$ is primitive and satisfies
\[
\tfrac34h_k\le h_{k+1}<\frac{(n-1-k)^2}{4(T-\delta_{U_k}^2)},
\]
together with Lemma~\ref{lem:kz}(b);
\item every $U_k$ has $v\in U_k^\perp$, so $U_k^\perp$ contains a vector with no zero coordinate;
\item if $k^\ast=\infty$, then $v$ lies in $\Z^n\cap U_{n-2}^\perp$, is primitive, and satisfies the bound
of Proposition~\ref{prop:leaf}.
\end{enumerate}
\end{theorem}

\begin{proof}
(i) is Proposition~\ref{prop:step} with $U=0$, together with two facts: shortest vectors are primitive, and
$\langle w_1,v\rangle=0$ with $v>0$. (ii) follows from Proposition~\ref{prop:step},
Lemma~\ref{lem:vfree}(b),(c) and Lemma~\ref{lem:kz}, since $h_{k+1}=\lambda_1(P_{U_k^\perp}L)^2$.
(iii) is immediate. (iv) is Proposition~\ref{prop:leaf}.
\end{proof}

The second alternative in the theorem cannot be ignored. But it is not a weakness either: it points to
directions that are genuinely hard to block.

\begin{proposition}[Unresolved nodes are witnesses]\label{prop:witness}
Let $U$ be a rational subspace with $\dim U\le n-2$ such that $U^\perp$ contains an integer vector with no
zero coordinate, and suppose $\delta_U^2\ge T$. Then there are infinitely many pairwise non-proportional
primitive $v\in\Z^n_{>0}$ with $d_n(v)\ge\delta_U^2\ge T$.
\end{proposition}

\begin{proof}
The integer vectors of $U^\perp$ with no zero coordinate form the complement of finitely many proper
subspaces in a lattice of rank $\ge2$. So there are infinitely many non-proportional primitive ones. For
each such $z$, Lemmas~\ref{lem:lower} and \ref{lem:sym}(a) give $d_n(|z|)=d_n(z)\ge\delta_U^2$.
\end{proof}

\begin{corollary}[Certificate]\label{cor:cert}
Consider a search that, up to a group preserving $d_n$ and all quantities in
Theorem~\ref{thm:complete}, does the following:
\begin{itemize}
\item it generates every root allowed by (i);
\item at every node with $\delta^2<T$, it generates every extension allowed by (ii) and (iii);
\item it \emph{records} every node with $\delta^2\ge T$ as unresolved, and does not expand it;
\item at every leaf it enumerates every vector allowed by (iv).
\end{itemize}
If the search records no unresolved node and every enumerated leaf vector $v$ satisfies $d_n(|v|)\le T$,
then $\Delta_n\le T$. If the leaf vectors with $d_n(|v|)>T'$ are exactly the permutations of $v^\ast$, then
every $v$ with $d_n(v)>T'$ is a permutation of $v^\ast$.
\end{corollary}

\section{The certified search}\label{sec:algorithm}

Corollary~\ref{cor:cert} says what a search must do to prove $\Delta_n\le T$. This section describes how we
carry it out. It covers two practical points: how to exploit the symmetry of $d_n$ without losing
positivity information, and how to make every decision in exact arithmetic.

\subsection{Symmetry reduction}

Let $G$ be the group of signed permutation matrices. It preserves $\Z^n$, the Euclidean metric and
$c+\Z^n$ (Lemma~\ref{lem:sym}). Hence it preserves offsets, saturations, quotient lattices, KZ conditions and
all bounds in Theorem~\ref{thm:complete}. It does not preserve positivity. We therefore replace
condition~(iii) by the $G$-invariant condition below, and take absolute values at the leaves.

\begin{lemma}[Full-support criterion]\label{lem:support}
Let $U$ be a rational subspace. The following are equivalent: $U^\perp$ contains an integer vector with no
zero coordinate; no standard basis vector $e_i$ lies in $U$; every diagonal entry of $P_U$ is $<1$. If
$r_1,\dots,r_k$ are pairwise orthogonal vectors spanning $U$, then $(P_U)_{ii}=\sum_jr_{j,i}^2/\|r_j\|^2$.
\end{lemma}

\begin{proof}
The coordinate $x_i$ vanishes identically on $U^\perp$ if and only if $e_i\in U$. A rational vector space is
not a finite union of proper subspaces, so some rational point of $U^\perp$ avoids every $\{x_i=0\}$.
Moreover, $(P_U)_{ii}=\|P_Ue_i\|^2$, which equals $1$ exactly when $e_i\in U$.
\end{proof}

\begin{lemma}[Complete orbit key]\label{lem:key}
Given a flag $U_1\subset\dots\subset U_k$ of rational subspaces, let $r_j$ be the primitive integer vector
spanning $U_j\cap U_{j-1}^\perp$; it is unique up to sign. Let $R$ be the $k\times n$ matrix with rows $r_j$.
Define $\kappa(R)$ as follows: over all row-sign patterns, normalize each column up to sign, sort the columns
lexicographically, and take the lexicographically least result. Then two flags are equivalent under $G$ if and
only if their keys $\kappa$ coincide.
\end{lemma}

\begin{proof}
An element $S\in G$ maps the flag $(U_j)$ onto $(U_j')$ if and only if $Sr_j=\pm r_j'$ for every $j$, because
$S$ is orthogonal and preserves primitivity. Such an $S$ acts on $R$ by permuting columns and changing their
signs. So two flags are equivalent exactly when their matrices lie in one orbit of the group generated by
row signs, column signs and column permutations. The key $\kappa$ is the least element of that orbit
under a fixed total order, so it is a complete invariant.
\end{proof}

\subsection{Algorithm}

Fix $n$ and $T$, and put $W=\lceil(n-1)^2/(4T)\rceil-1$. States are KZ-flag prefixes $(w_1,\dots,w_k)$, each
a basis of $\Z^n\cap U_k$, stored once per $G$-orbit using the key of Lemma~\ref{lem:key}.

\begin{enumerate}[label=\arabic*.]
\item \emph{Roots.} Every primitive $w$ with nonnegative sorted entries, at least two of them nonzero, and
$\|w\|^2\le W$. This covers condition (i) up to $G$.
\item \emph{Expansion of a state at depth $k<n-2$.}
  \begin{enumerate}[label=(\alph*)]
  \item Compute $\delta_{U_k}^2$ exactly. If $\delta_{U_k}^2\ge T$, record the state as \emph{unresolved} and stop.
  \item Put $\Lambda=(n-1-k)^2/(4(T-\delta_{U_k}^2))$. If $\Lambda\le\frac34h_k$, the state has no
  admissible extension (Lemma~\ref{lem:kz}(a)), so discard it.
  \item Using a Smith-form basis of $\Z^n/(\Z^n\cap U_k)$, enumerate every primitive class $\bar w$, up to
  sign, with $h=\|P_{U_k^\perp}w\|^2<\Lambda$. Keep those with $h\ge\frac34h_k$ that satisfy
  Lemma~\ref{lem:kz}(b); this is tested by exact closest-vector computations. Also require that
  $\spn(U_k,w)$ satisfy Lemma~\ref{lem:support}.
  \item For a child at depth $k+1<n-2$, compute its offset. If it is $\ge T$, record it as unresolved.
  If its next bound satisfies $(n-2-k)^2/(4(T-\delta^2))\le\frac34h$, discard it (Lemma~\ref{lem:kz}(a)
  one level ahead).
  \end{enumerate}
\item \emph{Leaves} (depth $n-2$). If $\delta^2\ge T$, record the leaf as unresolved. Otherwise enumerate
every primitive $z\in\Z^n\cap U_{n-2}^\perp$ with $\|z\|^2<\covol^2/(4(T-\delta^2))$ and no zero coordinate.
Evaluate $d_n(|z|)$ exactly with Lemma~\ref{lem:exact}.
\end{enumerate}

Every rejection above is a necessary condition from Theorem~\ref{thm:complete} or Lemma~\ref{lem:kz},
applied to the $G$-image of the KZ flag of a hypothetical counterexample. That image is again a KZ flag,
of $L_{Sv}$, and the leaf contains $Sv$, with $|Sv|$ a permutation of $v$. Hence the search satisfies the
hypotheses of Corollary~\ref{cor:cert}.

\subsection{Exactness}

All quantities are integers or rationals.
\begin{itemize}
\item Saturations, quotient bases and integer kernels come from an exact Smith normal form, whose
transformation matrices are checked.
\item Gram--Schmidt data, projectors, offsets and closest-vector distances are computed in $\Q$. The
closest-vector routine includes the component of the target orthogonal to the lattice.
\item Enumeration inside the ellipsoids $\{h<\Lambda\}$ is Fincke--Pohst/Schnorr--Euchner style
\cite{FP85,SE94}. Integer ranges are widened by a safety margin, and every point is accepted or rejected
by an exact rational comparison. The strict inequalities of Propositions~\ref{prop:step} and
\ref{prop:leaf} are respected.
\item LLL reduction \cite{LLL82} is used only to speed up enumeration. It does not affect which points are
found.
\end{itemize}

\section{Computations}\label{sec:computations}

We report three kinds of runs. First, runs at the classical values for $n\le5$ serve as a consistency check
against the literature. Second, the $n=6$ runs establish Theorem~A, using two independent programs and a
second run just below the threshold for uniqueness. Third, the $n=7$ run establishes Theorem~B.

The search of Section~\ref{sec:algorithm} is implemented in \texttt{m\_cvp\_global\_signs.py}, on top of the
exact routines in \texttt{m\_cvp\_proofsafe.py}. We also used an independent search mode of
\texttt{m\_cvp\_proofsafe.py}. It stores rational subspaces up to unsigned permutations, tests positivity
exactly through a Gordan--Stiemke alternative, and recomputes all admissible KZ flags of each subspace. Both
programs use only the Python standard library. Their SHA-256 digests are:
\begin{center}\scriptsize
\begin{tabular}{@{}ll@{}}
\texttt{m\_cvp\_proofsafe.py} & \texttt{d31f83f32912b1730fe965d5815ce92ed4e736b2e9788b63a5eb9ca003dbdffa}\\
\texttt{m\_cvp\_global\_signs.py} & \texttt{aff71c9072523bd44b0145b2cd4e590709702cb93a089d50de0d964b40f3fcce}
\end{tabular}
\end{center}

\subsection{Small dimensions}\label{subsec:small}

For $n=2$ no search is needed: $\|v\|^2=\|w_1\|^2\le4$ forces $v=(1,1)$, and $d_2(1,1)=0$.
Table~\ref{tab:small} records certified runs at the classical values.

\begin{table}[ht]
\caption{Certified runs for $n\le5$. Every run was verified and recorded no unresolved node.}\label{tab:small}
\centering\small
\begin{tabular}{llllr}
\toprule
$n$ & $T$ & program & states by depth & leaf candidates\\
\midrule
3 & $3/28$ & subspace & $7$ & $9$\\
4 & $11/60$ & subspace & $22,\,62$ & $57$\\
5 & $41/168$ & subspace & $61,\,1520,\,1365$ & $960$\\
5 & $41/168$ & signed & $24,\,79,\,618$ & $960$\\
\bottomrule
\end{tabular}
\end{table}

\subsection{Dimension six}

The runs for $n=6$ are listed in Table~\ref{tab:n6}.

\begin{table}[ht]
\caption{Certified runs for $n=6$, each on two processes. ${}^\ast$Subspace mode recomputing all admissible KZ flags.}\label{tab:n6}
\centering\small
\begin{tabular}{lllrr}
\toprule
$T$ & program & states by depth $1,2,3,4$ & leaves & time (s)\\
\midrule
$107/312$ & subspace & $118,\ 15185,\ 126235,\ 12112$ & $6168$ & $2107$\\
$107/312$ & signed & $33,\ 555,\ 1096,\ 6402$ & $6168$ & $222$\\
$107/312-\varepsilon$ & subspace$^\ast$ & $118,\ 4323,\ 7085,\ 25805$ & $6168$ & $1504$\\
\bottomrule
\end{tabular}
\end{table}

Here $\varepsilon=313/312{,}000{,}000$. Every run recorded $0$ unresolved nodes. The signed run's frontier
and candidate digests are \texttt{5d489f7e7920\ldots} and \texttt{aa284d9e4adb\ldots}.

\begin{proof}[Proof of Theorem A]
By Lemma~\ref{lem:exact}, $d_6(12,9,7,5,3,2)=107/312$. The runs at $T=107/312$ record no unresolved node,
and every leaf candidate has $d_6\le T$. Corollary~\ref{cor:cert} therefore gives $\Delta_6\le107/312$, and
hence equality. In the run at $T'=107/312-\varepsilon$ there is again no unresolved node, and the only leaf
candidate with $d_6>T'$ is $(12,9,7,5,3,2)$. By Corollary~\ref{cor:cert}, every $v$ with $d_6(v)=\Delta_6$ is a
permutation of it.
\end{proof}

\subsection{Dimension seven}

All $n=7$ runs were carried out on an Apple M1 Pro (8 cores: 6 performance and 2 efficiency; 16\,GB of
memory) under macOS~15.7.9 with Python~3.9.6, using $8$ worker processes. The runs for $n\le6$ used a
2-core Linux container with Python~3.11. The first $n=7$ run used $T=73/180$ and $W=22$; its total wall
time was $96{,}493$\,s (Table~\ref{tab:n7}).

\begin{table}[ht]
\caption{First run for $n=7$ at $T=73/180$. The last column gives the prefix of the frontier SHA-256 digest.}\label{tab:n7}
\centering\small
\begin{tabular}{lrrrrl}
\toprule
depth & in & out & enumerated & time (s) & digest\\
\midrule
$1\to2$ & $58$ & $6{,}069$ & $1{,}931{,}321$ & $82$ & \texttt{edeabb6c865b}\\
$2\to3$ & $6{,}069$ & $154{,}363$ & $59{,}218{,}407$ & $7{,}230$ & \texttt{1595e1fffac8}\\
$3\to4$ & $154{,}363$ & $66{,}272$ & $291{,}079{,}485$ & $81{,}491$ & \texttt{ff18fb31fab1}\\
$4\to5$ & $66{,}272$ & $392{,}264$ & $16{,}055{,}230$ & $2{,}512$ & \texttt{297b5e3deb4f}\\
leaves & $392{,}264$ & -- & -- & $4{,}961$ & \\
\bottomrule
\end{tabular}
\end{table}

The leaf stage merged the $392{,}264$ leaf flags into $328{,}883$ distinct subspaces, up to unsigned
permutations. It enumerated $9{,}543{,}063$ kernel vectors, giving $403{,}545$ distinct candidates (digest
\texttt{429dc9cd4ffb\ldots}), and evaluated all of them exactly. No unresolved node was recorded, and no
candidate has $d_7>73/180$. The depth-$2$ frontier was regenerated independently on a second machine, and
its digest agreed.

\emph{Independent reproduction.} The complete $n=7$ computation at $T=73/180$ was then repeated with a
separate driver. It processes each layer in resumable chunks and merges them by the orbit key of
Lemma~\ref{lem:key}. Its leaf stage evaluates all $392{,}264$ leaf flags individually, with no merging of
equal subspaces, which gives $11{,}400{,}280$ kernel vectors instead of $9{,}543{,}063$. Candidates are
first screened at exact rational sample points $t=a/m$; a value $\le T$ there proves $d_7\le T$. Every
survivor is then evaluated exactly. This second run took about $25$ hours of wall time with $8$ workers.
All four frontier digests, all pruning counts, the $403{,}545$ unique candidates and their digest
\texttt{429dc9cd4ffb\ldots} agree exactly with the first run. It again recorded no unresolved node and no
candidate with $d_7>73/180$.

\emph{Uniqueness run.} For uniqueness the search was repeated with the second driver at
\[
T'=\frac{73}{180}-10^{-6}=\frac{3649991}{9000000}.
\]
The root bound is unchanged, since $W=22$ for every $\varepsilon<73/180-9/23\approx0.0143$. The run took
about $24$ hours of wall time with $8$ workers (Table~\ref{tab:n7u}).

\begin{table}[ht]
\caption{Uniqueness run for $n=7$ at $T'=3649991/9000000$. The last column gives the prefix of the frontier SHA-256 digest.}\label{tab:n7u}
\centering\small
\begin{tabular}{lrrrrl}
\toprule
depth & in & out & enumerated & time (s) & digest\\
\midrule
$1\to2$ & $58$ & $6{,}069$ & $1{,}931{,}321$ & $77$ & \texttt{edeabb6c865b}\\
$2\to3$ & $6{,}069$ & $154{,}501$ & $59{,}218{,}855$ & $6{,}544$ & \texttt{1d0f4418fcbc}\\
$3\to4$ & $154{,}501$ & $66{,}404$ & $291{,}352{,}222$ & $71{,}840$ & \texttt{c47eaea81008}\\
$4\to5$ & $66{,}404$ & $392{,}690$ & $16{,}092{,}668$ & $1{,}344$ & \texttt{62c3876be0a7}\\
leaves & $392{,}690$ & -- & -- & $5{,}589$\,\textsuperscript{$\dagger$} & \\
\bottomrule
\end{tabular}

\smallskip
{\footnotesize $^\dagger$\,Sum of per-chunk wall times.}
\end{table}

The depth-$2$ frontier coincides with the one at $T=73/180$, and the later frontiers are slightly larger,
as expected for a lower threshold. The leaf stage enumerated $11{,}414{,}179$ kernel vectors, giving
$404{,}230$ distinct candidates (digest \texttt{0b93524f8f6a\ldots}). No unresolved node was recorded.
Exactly one candidate has $d_7>T'$, namely $(11,10,8,7,4,3,1)$, with $d_7=73/180$.

\begin{proof}[Proof of Theorem B]
By Lemma~\ref{lem:exact}, $d_7(11,10,8,7,4,3,1)=73/180$. The runs at $T=73/180$ record no unresolved node,
and every leaf candidate has $d_7\le T$. Corollary~\ref{cor:cert} therefore gives $\Delta_7\le73/180$, and
hence equality. In the run at $T'=3649991/9000000$ there is again no unresolved node, and the only leaf
candidate with $d_7>T'$ is $(11,10,8,7,4,3,1)$. By Corollary~\ref{cor:cert}, every $v$ with
$d_7(v)=\Delta_7$ is a permutation of it.
\end{proof}

\begin{remark}\label{rem:unique7}
Both uniqueness runs show more than uniqueness. Every positive primitive $v$ that is not a permutation of
the extremal vector satisfies $d_6(v)\le107/312-10^{-6}$, respectively $d_7(v)\le73/180-10^{-6}$. So the
maximum is isolated by a gap of at least $10^{-6}$.
\end{remark}

\subsection{Validation}

\begin{itemize}
\item Built-in regression tests cover Smith forms, integer kernels, the orthogonal residual in closest-vector
computations, invariance of the orbit keys, and KZ flags of ranks $3$ and $4$.
\item The evaluator of Lemma~\ref{lem:exact} agrees with a separate implementation, which clips each cell's
minimizer to the cell, on random vectors.
\item Offsets agree with brute-force minimization over boxes.
\item The inequality of Lemma~\ref{lem:lower} held on $3{,}632$ random pairs $(v,U)$. The bounds of
Propositions~\ref{prop:step} and \ref{prop:leaf} held on $1{,}193$ and $690$ random instances.
\item For $n=5,6$ the two programs, which use different symmetry reductions, positivity tests and state
representations, produced identical sets of leaf candidates.
\item A direct search over all positive $v\in\Z^6$ with $\max_iv_i\le24$ found no vector other than
$(12,9,7,5,3,2)$ with $d_6(v)\ge107/312-10^{-6}$.
\end{itemize}

\section{Remarks}\label{sec:remarks}

\begin{remark}[No dimensional barrier]
Every step of the search is finite, and Proposition~\ref{prop:witness} converts any unresolved node into
vectors with $d_n\ge T$. The method is therefore complete in every dimension, up to a single exception. Suppose
$T=\Delta_n$ and some rational subspace $U$ with $\dim U\le n-2$, whose orthogonal complement contains a
vector with no zero coordinate, had $\delta_U^2=T$ exactly. Then there would be infinitely many
maximizers, and the search could not close. Beyond that, the obstruction is
cost. In the $n=7$ run, the transition $3\to4$ dominates: it examined $2.9\times10^8$ quotient vectors,
of which $81\%$ failed the Gauss height condition and $14\%$ the KZ test.
\end{remark}

\begin{remark}[The sequence $\Delta_n$]
By Lemma~\ref{lem:sym}(b), the sequence $(\Delta_n)$ is nondecreasing. The known values
$\frac1{20},\frac3{28},\frac{11}{60},\frac{41}{168},\frac{107}{312},\frac{73}{180}$ are strictly increasing,
and $\Delta_n/(n/12)$ increases from $0.30$ to about $0.70$ over $2\le n\le7$
(see Table~\ref{tab:values}). Its behaviour for large $n$ is not addressed here.
\end{remark}

\begin{remark}[Beyond $n=7$]
Nothing in Sections~\ref{sec:lattices}--\ref{sec:algorithm} depends on the dimension. The step bounds weaken
as $n$ grows, because the factor $(n-1-k)^2$ increases while the relevant thresholds grow only slowly. The
search tree therefore grows quickly. For $n=7$, the depth-$3\to4$ layer alone accounted for
$85\%$ of the running time. A candidate threshold for $n=8$ and an estimate of the cost of the
corresponding search are natural next steps.
\end{remark}

\section*{Data availability}

The two programs, all checkpoints (frontiers at every depth), leaf-candidate files and run logs are available
at \url{https://doi.org/10.5281/zenodo.23130200} and \url{https://github.com/vineet-kahlon/View-Obstruction-Problem}. Each frontier can be checked independently against the digests above.

\section*{Acknowledgments}

The author used a Gemini-based AI assistant for help with implementing and
optimizing the search code and for language editing. The mathematical results,
proofs and the design of the search are the author's own.
The computations were
reproduced independently as described in Section~\ref{sec:computations}.

\appendix
\section{Full frontier digests for \texorpdfstring{$n=7$}{n=7}}

{\scriptsize\ttfamily
\begin{tabular}{ll}
depth 2 & edeabb6c865b89845a8b3c558287abd8aacc54ff9616c05a30e17cb0a32dc96c\\
depth 3 & 1595e1fffac8ba18d07677fd67f83b4412830a54815b911dee4870c43cfe4e88\\
depth 4 & ff18fb31fab1908b24632443dae9128658d742217003b70b45276e3544a2cf16\\
depth 5 & 297b5e3deb4faf9366d6a1f2f73ef026a7166496a280c2050db673ae07e1e819\\
candidates & 429dc9cd4ffbcb148f3b316d8f4fcbb8eb72274153d3773d373ee3910caa7e07\\[3pt]
\multicolumn{2}{l}{\normalfont\small Uniqueness run, $T'=3649991/9000000$:}\\
depth 2 & edeabb6c865b89845a8b3c558287abd8aacc54ff9616c05a30e17cb0a32dc96c\\
depth 3 & 1d0f4418fcbcce612afeea917c4dba6a28c9cdfbafecea50b365a9df9143d41f\\
depth 4 & c47eaea81008e47e2886416d34955dc3957d63a827889ce867a905b421535af9\\
depth 5 & 62c3876be0a710fcebeee28fd64421202b0c50ab08d49c0cea000e0b56a72932\\
candidates & 0b93524f8f6a5bc246dc69a93b25682560049e1012fc4617b00b5d6ea5ab7d27
\end{tabular}}


\end{document}